\documentclass[12pt]{iopart}
\usepackage{amsfonts}
\usepackage[english]{babel}
\usepackage{euscript}
\usepackage{bbm}

\begin{document}

\title[Deformed cohomologies of pseudo-groups and coverings of differential equations]%
{Deformed cohomologies of symmetry pseudo-groups and coverings of differential equations}

\author{Oleg I. Morozov}

\address{Faculty of Applied
  Mathematics, AGH University of Science and Technology,
  \\
  Al. Mickiewicza 30,
  Krak\'ow 30-059, Poland
\\
morozov{\symbol{64}}agh.edu.pl}

\begin{abstract}
We discuss a relation between deformed cohomologies  of symmetry pseu\-do-groups and coverings of differential equations. Examples include
the potential Khokhlov--Zabolotskaya equation and the Boyer--Finley equa\-tion.
\end{abstract}



\section{Introduction}

Deformed (or exotic) cohomologies were introduced  in works of S.P. Novikov, \cite{Novikov1986,Novikov2002,Novikov2005}, as a tool in an analogue of the Morse theory for smooth multi-valued fun\-c\-ti\-ons. Then they were applied to different problems of symplectic geometry and algebraic topology, \cite{Alaniya1997,Alaniya1999,Milionshchikov2002,Milionshchikov2005,Pazhitnov1991}. The objective of the present paper is to establish a relation between the deformed cohomologies of symmetry pseudo-groups of partial differential equations and their coverings.

Coverings (or Wahlquist--Estabrook prolongation structures, \cite{WE}, or zero-cur\-va\-tu\-re representations, \cite{ZakharovShabat}, or integrable extensions, \cite{BryantGriffiths}, etc.) are of great im\-por\-tan\-ce in geometry of {\sc pde}s. The theory of coverings is
a natural  framework for dealing with nonlocal symmetries and nonlocal conservation laws, inverse scattering con\-struc\-ti\-ons for so\-li\-ton equations, B\"acklund transformations, recursion operators, and de\-for\-ma\-ti\-ons of non\-li\-ne\-ar {\sc pde}s, \cite{KV84,KV89,KV99}.
A number of techniques has been devised to handle the problem of recognizing whether a given differential equation has a covering,
\cite{WE,Morris1976,Morris1979,Estabrook,Zakharov82,DoddFordy,Tondo,Hoenselaers,Marvan1992,Sakovich,Marvan1997,Marvan2002}.
In \cite{Kuzmina}, examples of coverings of {\sc pde}s with three independent variables were found by means of \'Elie Cartan's method of equivalence, \cite{Cartan1,Cartan2,Cartan4,Gardner,Kamran,Olver95}.
This idea was developed in \cite{Morozov2007,Morozov2008,Morozov2010}. In \cite{Morozov2009} we propose an approach to the covering problem based on the technique of contact in\-te\-gra\-ble extensions ({\sc cie}s) of the structure equations of the symmetry pseudo-groups, which is a generalization of the definition of integrable extension from \cite[\S6]{BryantGriffiths} for the case of more than two independent variables. Then in \cite{Morozov2010b,Morozov2012,Morozov2014,Morozov2014b} the method of {\sc cie}s was applied to finding of coverings, B\"acklund transformations and recursion operators for a number of {\sc pde}s.

From the definition of deformed cohomology it follows that each non-trivial de\-for\-med 2-co\-cy\-cle of the symmetry pseudo-group of a {\sc pde} provides an integrable ex\-ten\-si\-on of this pseudo-group. Then a covering for the {\sc pde} may be obtained via integration of the extension equation in accordance with Cartan's theorem.

In this paper we consider two equations:  the potential Khokhlov--Zabolotskaya equation (or Lin--Reissner--Tsien equation), \cite{LinReissnerTsien,KhZ},
\begin{equation}
u_{yy} = u_{tx} + u_x\,u_{xx},
\label{pKhZ}
\end{equation}
and the Boyer--Finley equation, \cite{BoyerFinley},
\begin{equation}
u_{tx} = \e^{u_y} u_{yy}.
\label{BF}
\end{equation}
We show that symmetry pseudo-groups of both equations have non-trivial deformed second cohomologies.
The integrable extensions that correspond to cocycles from these cohomology groups define known coverings of equations
(\ref{pKhZ}) and (\ref{BF}).

\section{Preliminaries}\label{Preliminaries_section}

\subsection{Coverings of PDEs}

All considerations in this paper are local. The presentation in this subsection closely follows to
\cite{KrasilshchikVerbovetsky2011,KrasilshchikVerbovetskyVitolo2012}.
Let $\pi \colon \mathbb{R}^n \times \mathbb{R}^m \rightarrow \mathbb{R}^n$,
$\pi \colon (x^1, \dots, x^n, u^1, \dots, u^m) \mapsto (x^1, \dots, x^n)$ be a trivial bundle, and $J^\infty(\pi)$
be the bundle of its jets of the infinite order. The local coordinates on $J^\infty(\pi)$ are $(x^i,u^\alpha,u^\alpha_I)$,
where $I=(i_1, \dots, i_n)$ is a multi-index, and for every local section
$f \colon \mathbb{R}^n \rightarrow \mathbb{R}^n \times \mathbb{R}^m$ of $\pi$ the corresponding infinite jet $j_\infty(f)$ is
a section $j_\infty(f) \colon \mathbb{R}^n \rightarrow J^\infty(\pi)$ such that
$u^\alpha_I(j_\infty(f))
=\displaystyle{\frac{\partial ^{\#I} f^\alpha}{\partial x^I}}
=\displaystyle{\frac{\partial ^{i_1+\dots+i_n} f^\alpha}{(\partial x^1)^{i_1}\dots (\partial x^n)^{i_n}}}$.
We put $u^\alpha = u^\alpha_{(0,\dots,0)}$. Also, in the case of $n=3$, $m=1$ we denote $x^1 = t$, $x^2= x$,
$x^3= y$, and $u^1_{(i,j,k)}=u_{{t \dots t}{x \dots x}{y \dots y}}$
with $i$  times $t$, $j$  times $x$, and $k$  times $y$.

The vector fields
\[
D_{x^k} = \frac{\partial}{\partial x^k} + \sum \limits_{\# I \ge 0} \sum \limits_{\alpha = 1}^m
u^\alpha_{I+1_{k}}\,\frac{\partial}{\partial u^\alpha_I},
\qquad k \in \{1,\dots,n\},
\]
$(i_1,\dots, i_k,\dots, i_n)+1_k = (i_1,\dots, i_k+1,\dots, i_n)$,  are called {\it total derivatives}.  They commute everywhere on
$J^\infty(\pi)$:  $[D_{x^i}, D_{x^j}] = 0$.

\vskip 10 pt

%

A system of {\sc pde}s $F_r(x^i,u^\alpha_I) = 0$, $\# I \le s$, $r \in \{1,\dots, R\}$, of the order
$s \ge 1$ with $R \ge 1$ defines the submanifold
$\EuScript{E} = \{(x^i,u^\alpha_I) \in J^\infty(\pi) \,\,\vert\,\, D_K(F_r(x^i,u^\alpha_I)) = 0, \,\, \# K \ge 0\}$
in $J^\infty(\pi)$.

\vskip 10 pt


Denote $\mathcal{W} = \mathbb{R}^\infty$ with  coordinates $w^s$, $s \in  \mathbb{N} \cup \{0\}$. Locally,
an (infinite-di\-men\-si\-o\-nal)  {\it differential covering} of $\EuScript{E}$ is a trivial bundle $\tau \colon J^\infty(\pi) \times \mathcal{W} \rightarrow J^\infty(\pi)$
equipped with the {\it extended total derivatives}
\begin{equation}
\tilde{D}_{x^k} = D_{x^k} + \sum \limits_{ s =0}^\infty
T^s_k(x^i,u^\alpha_I,w^j)\,\frac{\partial }{\partial w^s}
\label{extended_derivatives}
\end{equation}
such that $[\tilde{D}_{x^i}, \tilde{D}_{x^j}]=0$ for all $i \not = j$ whenever $(x^i,u^\alpha_I) \in \EuScript{E}$. We define
the partial derivatives of $w^s$ by  $w^s_{x^k} =  \tilde{D}_{x^k}(w^s)$.  This yields the system of
{\it covering equations}
\[
w^s_{x^k} = T^s_k(x^i,u^\alpha_I,w^j).
\]
This over-determined system of {\sc pde}s is compatible whenever $(x^i,u^\alpha_I) \in \EuScript{E}$.

\vskip 5 pt

Dually the covering with extended total derivatives (\ref{extended_derivatives}) is defined by the integrable ideal of
the {\it Wahlquist--Estabrook forms}
\[
dw^s - T^s_k(x^i,u^\alpha_I,w^j)\,dx^k.
\]

\subsection{Cartan's structure theory of Lie pseudo-groups}

Let $M$ be a manifold of dimension $n$. A {\it local diffeomorphism} on $M$ is a diffeomorphism
$\Phi \colon \EuScript{U} \rightarrow \hat{\EuScript{U}}$ of two open subsets of $M$. A {\it pseudo-group}
$\mathfrak{G}$ on $M$ is a collection of local dif\-feo\-mor\-phisms of $M$, which is closed under composition
whenever the latter is defined, contains an identity and is closed under inversion. A {\it Lie pseudo-group} is a pseudo-group
whose diffeomorphisms are local analytic solutions of an involutive system of partial differential equations called
{\it defining system}.

\'Elie Cartan's approach to Lie pseudo-groups is based on a possibility to characterize transformations from
a pseudo-group in terms of a set of invariant differential 1-forms called {\it Maurer--Cartan} ({\sc mc}) {\it forms}.
In a general case, {\sc mc} forms $\omega^1$, ... , $\omega^m$ of an infinite-dimensional Lie pseudo-group
$\mathfrak{G}$ are defined on a direct product $M \times \tilde{M} \times G$, where $\tilde{M}$ is the coordinate
space of parameters of prolongation, \cite[Ch. 12]{Olver95}, $G$ is a finite-dimensional Lie group, and
$m = \dim \,M +\dim \,\tilde{M}$. The forms $\omega^i$ are independent and include differentials of coordinates on
$M \times \tilde{M}$ only, while their coefficients depend also on coordinates of $G$. These  forms characterize the
pseudo-group $\mathfrak{G}$ in the following sense: a local diffeomorphism
$\Phi \colon \EuScript{U} \rightarrow \hat{\EuScript{U}}$ on $M$ belongs to $\mathfrak{G}$ whenever there exists a
local diffeomorphism $\Psi \colon \EuScript{W} \rightarrow \hat{\EuScript{W}}$ on $M \times \tilde{M}\times G$ such
that $\upsilon \circ \Psi = \Phi \circ \upsilon$ for the projection
$\upsilon \colon M \times \tilde{M} \times G \rightarrow M$ and the forms $\omega^j$ are invariant w.r.t. $\Psi$, that
is,
\begin{equation}
\Psi^{*} \left(\omega^i\vert {}_{\hat{\EuScript{W}}} \right)
= \omega^i\vert {}_{\EuScript{W}}.
\label{Phi_omega}
\end{equation}
Expressions for $d \omega^i$ in terms of $\omega^j$ give Cartan's
{\it structure equations} of $\mathfrak{G}$:
\begin{equation}
d \omega^i = A_{\gamma j}^i \,\pi^\gamma \wedge \omega^j + \case{1}{2}\,B_{jk}^i\,\omega^j \wedge \omega^k,
\qquad B_{jk}^i = - B_{kj}^i.
\label{SE}
\end{equation}
The forms $\pi^\gamma$, $\gamma \in \{1,...,\dim \, G\}$, are linear combinations of {\sc mc} forms of the Lie group
$G$ and the forms $\omega^i$. The coefficients $A_{\gamma j}^i$ and $B_{jk}^i$ are either constants or functions of
a set of invariants $U^\kappa \colon M \rightarrow \mathbb{R}$, $\kappa \in \{1,...,l\}$,
$l < \dim\, M$, of the pseudo-group $\mathfrak{G}$, so
$\Phi^{*} \left(U^{\kappa}\vert {}_{\hat{\EuScript{U}}} \right) = U^{\kappa}\vert {}_{\EuScript{U}}$
for every $\Phi \in \mathfrak{G}$. In the latter case, the differentials of $U^\kappa$ are invariant 1-forms, so
they are linear combinations of the forms $\omega^j$,
\begin{equation}
d U^\kappa = C_j^\kappa\,\omega^j,
\label{dUs}
\end{equation}
where the coefficients $C_j^\kappa$ depend on the invariants $U^1$, ..., $U^l$ only.

Equations (\ref{SE}) must be compatible in the following sense: we have
\begin{equation}
d(d \omega^i) = 0 = d \left(A_{\gamma j}^i \,\pi^\gamma \wedge \omega^j + \case{1}{2}\,B_{jk}^i\,\omega^j \wedge \omega^k \right),
\label{compatibility_conditions_SE}
\end{equation}
therefore there must exist expressions
\begin{equation}
d \pi^\gamma = W_{\lambda j}^\gamma\, \chi^\lambda \wedge \omega^j
+ X_{\beta \epsilon}^\gamma\,\pi^\beta\wedge\pi^\epsilon
+Y_{\beta j}^\gamma\,\pi^\beta \wedge \omega^j
+Z_{jk}^\gamma\,\omega^j \wedge \omega^k
\label{prolonged_SE}
\end{equation}
with some additional 1-forms $\chi^\lambda$  such that the right-hand side of (\ref{compatibility_conditions_SE})
is identically equal to zero after substituting for (\ref{SE}), (\ref{dUs}), and (\ref{prolonged_SE}).
Also, from (\ref{dUs}) it follows that
the right-hand side of the equation
\begin{equation}
d(d U^\kappa) = 0 = d(C_j^\kappa\,\omega^j)
\label{compatibility_conditions_dUs}
\end{equation}
must be identically equal to zero after substituting for (\ref{SE}) and (\ref{dUs}).

The forms $\pi^\gamma$ are not invariant w.r.t. the pseudo-group $\mathfrak{G}$. Respectively, the structure equations (\ref{SE}) are not changing when replacing
$\pi^\gamma \mapsto \pi^\gamma + z^\gamma_j\,\omega^j$ for certain parametric coefficients $z^\gamma_j$. The dimension $r^{(1)}$ of the linear space of these coefficients satisfies the fol\-lo\-wing inequality
\begin{equation}
r^{(1)} \le n\,\dim \, G  -  \sum \limits_{k=1}^{n-1} (n-k)\,s_k,
\label{CT}
\end{equation}
where the {\it reduced characters} $s_k$ are defined by the formulas
\begin{eqnarray}
s_1 &=& \max \limits_{u_1\in \mathbb{R}^n}\, \mathrm{rank}\,\, \mathbb{A}_1(u_1),
\nonumber
\\
s_k &=& \max \limits_{u_1,...,u_k \in \mathbb{R}^n}\, \mathrm{rank}\,\, \mathbb{A}_k(u_1,...,u_k) -
\sum \limits_{j=1}^{k-1} s_j, \qquad k \in \{1, ... , n-1\},
\nonumber
\\
s_n &=& \mathrm{dim}\, G - \sum \limits_{j=1}^{n-1} s_j,
\nonumber
\end{eqnarray}
with the  matrices $\mathbb{A}_k$ inductively defined by
\[
\fl
\mathbb{A}_1(u_1) = \left(A^i_{\gamma j} \,u^j_1\right),
\qquad \mathbb{A}_l(u_1,...,u_l) = \left(
\begin{array}{c}
\mathbb{A}_{l-1}(u_1,...,u_{l-1})
\\
A^i_{\gamma j} \,u^j_l
\end{array}
\right),
\qquad l \in \{2, ... n-1\},
\]
see \cite[\S 5]{Cartan1}, \cite[Def. 11.4]{Olver95} for the full discussion. The system of forms $\omega^k$ is {\it involutive}  when both sides of (\ref{CT}) are equal, \cite[\S 6]{Cartan1}, \cite[Def. 11.7]{Olver95}.

Cartan's fundamental theorems, \cite[\S\S 16, 22--24]{Cartan1}, \cite{Cartan4},
\cite[\S\S 16, 19, 20, 25,26]{Vasilieva1972}, \cite[\S\S 14.1--14.3]{Stormark2000}, state that for a Lie
pseudo-group there exists a set of {\sc mc} forms whose structure equations satisfy the compatibility and
involutivity conditions; conversely, if equations (\ref{SE}), (\ref{dUs}) meet the compatibility conditions
(\ref{compatibility_conditions_SE}), (\ref{compatibility_conditions_dUs}) and the involutivity con\-di\-tion, then
there exists a collection of 1-forms $\omega^1$, ... , $\omega^m$ and functions $U^1$, ... , $U^l$ which satisfy
(\ref{SE}) and (\ref{dUs}).  Equations (\ref{Phi_omega}) then define local diffeomorphisms from a Lie pseudo-group.

\vskip 10 pt
\noindent
{\sc Example 1}.
Suppose $\EuScript{E}$ is a second-order differential equation in one dependent and $n$ independent variables. We consider $\EuScript{E}$ as a submanifold in $J^2(\pi)$ with $\pi \colon \mathbb{R}^n \times \mathbb{R} \rightarrow \mathbb{R}^n$.
Let ${\rm{Cont}}(\EuScript{E})$ be the group of contact symmetries for $\EuScript{E}$. It consists of all the contact transformations on $J^2(\pi)$ mapping $\EuScript{E}$ to itself.
%
%
The {\sc mc} forms of ${\rm{Cont}}(\EuScript{E})$ can be computed from  the {\sc mc} forms
of the pseudo-group all the contact transformations on $J^2(\pi)$
 algorithmically by means of Cartan's method of equivalence, \cite{Cartan1,Cartan2,Cartan4,Gardner,Kamran,Olver95}, see details and examples in \cite{FelsOlver,Morozov2002,Morozov2006}.

\subsection{Deformed cohomologies}

Let $\mathfrak{g}$ be a Lie algebra over $\mathbb{R}$ and $\rho \colon \mathfrak{g} \rightarrow \mathrm{End}(V)$ be its representation. Let $C^k(\mathfrak{g}, V) =\mathrm{Hom}(\Lambda^k(\mathfrak{g}), V)$, $k \ge 1$,  be the space of all $k$--linear skew-symmetric mappings from $\mathfrak{g}$ to $V$. Then a differential complex
\[
V=C^0(\mathfrak{g}, V) \stackrel{d}{\longrightarrow} C^1(\mathfrak{g}, V)
\stackrel{d}{\longrightarrow} \dots \stackrel{d}{\longrightarrow}
C^k(\mathfrak{g}, V) \stackrel{d}{\longrightarrow} C^{k+1}(\mathfrak{g}, V)
\stackrel{d}{\longrightarrow} \dots
\]
is defined by the formula
\begin{eqnarray*}
\fl
d \theta (X_1, ... , X_{k+1}) &=&
\sum\limits_{q=1}^{k+1}
(-1)^{q+1} \rho (X_q)\,(\theta (X_1, ... ,\hat{X}_q, ... ,  X_{k+1}))
\\
&&+\sum\limits_{1\le p < q \le k+1} (-1)^{p+q}
\theta ([X_p,X_q],X_1, ... ,\hat{X}_p, ... ,\hat{X}_q, ... ,  X_{k+1}).
\end{eqnarray*}
Cohomologies of the complex $(C^{*}(\mathfrak{g}, V), d)$ are referred to as {\it cohomologies of the Lie algebra} $\mathfrak{g}$
{\it with coefficents in the representation} $\rho$.
For the trivial representation $\rho_0 \colon \mathfrak{g} \rightarrow \mathbb{R}$,
$\rho_0 \colon X \mapsto 0$, cohomologies of the corresponding complex are called
{\it cohomologies with trivial coefficents} and denoted by $H^{*}(\mathfrak{g})$.

Consider a Lie algebra $\mathfrak{g}$ over $\mathbb{R}$ with non-trivial first cohomology group $H^1(\mathfrak{g})$
and take a closed 1-form $\omega$ on $\mathfrak{g}$. Then for any $\lambda \in \mathbb{R}$ define new {\it deformed differential} $d_{\lambda\omega} \colon C^k(\mathfrak{g},\mathbb{R}) \rightarrow C^{k+1}(\mathfrak{g},\mathbb{R})$ by
the formula
\[
d_{\lambda \omega} \theta = d \theta +\lambda \,\omega \wedge \theta.
\]
From  $d\omega = 0$ it follows that
\begin{equation}
d_{\lambda \omega} ^2=0.
\label{d_deformed_2}
\end{equation}
Cohomologies of the complex
\[
C^1(\mathfrak{g}, \mathbb{R})
\stackrel{d_{\lambda \omega}}{\longrightarrow}
\dots
\stackrel{d_{\lambda \omega}}{\longrightarrow}
C^k(\mathfrak{g}, \mathbb{R})
\stackrel{d_{\lambda \omega}}{\longrightarrow}
C^{k+1}(\mathfrak{g}, \mathbb{R})
\stackrel{d_{\lambda \omega}}{\longrightarrow} \dots
\]
are referred to as {\it deformed} (or {\it exotic}) {\it cohomologies} of $\mathfrak{g}$ and denoted by
$H^{*}_{\lambda\omega}(\mathfrak{g})$.

\vskip 5 pt

\noindent
{\sc Remark 1}.
Cohomologies $H^{*}_{\lambda\omega}(\mathfrak{g})$ coincide with cohomologies of $\mathfrak{g}$ with coefficients in
the one-dimensional representation $\rho_{\lambda\omega} \colon \mathfrak{g} \rightarrow \mathbb{R}$,
$\rho_{\lambda\omega} \colon X \mapsto \lambda\, \omega(X)$. In particular,
when $\lambda=0$, cohomologies $H^{*}_{\lambda\omega}(\mathfrak{g})$ coincide with $H^{*}(\mathfrak{g})$.

\vskip 5 pt

\noindent
{\sc Remark 2}.
In general, $d_{\lambda \omega} (\alpha \wedge \beta) \neq d_{\lambda \omega} (\alpha) \wedge \beta
+ (-1)^{\mathrm{deg}\, \alpha} \alpha \wedge d_{\lambda \omega} (\beta)$.

\vskip 5 pt
\noindent
{\sc example 2}.  Consider the system
\[
\fl
d\theta^1 = 0,
\quad
d\theta^2 = - \theta^1 \wedge \theta^2,
\quad
d\theta^3 =  \theta^1 \wedge \theta^3,
\quad
d\theta^4 =  2\,\theta^1 \wedge \theta^4,
\quad
d\theta^5 = \theta^2 \wedge \theta^3
\]
for 1-forms $\theta^1$, ... , $\theta^5$. This system is compatible, that is, applying $d$ to both sides of its equations and then substituting for the
equations themselves into the right-hand sides gives identities $0=0$. Therefore the system defines Lie algebra $\mathfrak{h}$ of vectors $X_1$, ... , $X_5$ such that $\theta^i(X_j)= \delta^i_j$. The forms $\theta^i$ are Maurer--Cartan forms of $\mathfrak{h}$.
Evidently, $H^1(\mathfrak{h}) = \mathbb{R}\,[\theta^1] =\mathbb{R}\,\theta^1$.
Then direct computations give:
\[
\fl
H^2_{\lambda\,\theta^1}(\mathfrak{h}) =
\left\{
\begin{array}{lllll}
\{0\} && \mathrm{for}&  \lambda \not \in \{-3, -2, -1, 1\},
\\
\mathbb{R}\,[\theta^3 \wedge \theta^4] && \mathrm{for}&  \lambda=-3,
\\
\mathbb{R}\,[\theta^1 \wedge \theta^4] && \mathrm{for}&  \lambda=-2,
\\
\mathbb{R}\,[\theta^1 \wedge \theta^3] \oplus \mathbb{R}\,[\theta^2 \wedge \theta^4]
\oplus \mathbb{R}\,[\theta^3 \wedge \theta^5]
&& \mathrm{for}&  \lambda=-1,
\\
\mathbb{R}\,[\theta^1 \wedge \theta^2] \oplus \mathbb{R}\,[\theta^2 \wedge \theta^5]
&& \mathrm{for}&  \lambda=1.
\end{array}
\right.
\]

\section{Symmetry pseudo-groups of the potential Khoklov-Zabolotskaya and the Boyer-Finley equations}

Using the procedures of \'Elie Cartan's method of equivalence 
 we find the Maurer--Cartan forms and their structure equations for the symmetry pseudo-groups of equations (\ref{pKhZ}) and  (\ref{BF}), see notation in  \cite{Morozov2006}.

\subsection{Potential Khoklov-Zabolotskaya equation}
The structure equations for the symmetry pseudo-group of equation (\ref{pKhZ}) read
\[
\fl
d\theta_0 =
\eta_1 \wedge \theta_0
+ \xi^1 \wedge \theta_1
+ \xi^2 \wedge \theta_2
+ \xi^3 \wedge \theta_3,
\]
\[
\fl
d\theta_1 =
 (\eta_1-\eta_2) \wedge \theta_1
-2\, \eta_3 \wedge \theta_3
- \theta_0 \wedge (\xi^2+\sigma_{22})
+\xi^1 \wedge \sigma_{11}
+\xi^2 \wedge \sigma_{12}
+\xi^3 \wedge \sigma_{13},
\]
\[
\fl
d\theta_2 =
 \case{1}{2}  \, (\eta_1 -\eta_2)\wedge \theta_2
+ \xi^1 \wedge \sigma_{12}
+ \xi^2 \wedge \sigma_{22}
+ \xi^3 \wedge \sigma_{23},
\]
\[
\fl
d\theta_3 =
 \case{3}{4}  \, (\eta_1-\eta_2) \wedge 3\,\theta_3
- \eta_3 \wedge \theta_2
+ \xi^1 \wedge \sigma_{13}
+ \xi^2 \wedge \sigma_{23}
+ \xi^3 \wedge \sigma_{12},
\]
\[
\fl
d\xi^1 =
 \eta_2 \wedge \xi^1,
\]
\[
\fl
d\xi^2 =
 \case{1}{2}  \, (\eta_1+\eta_2) \wedge \xi^2
+ \eta_3 \wedge \xi^3
- \theta_2 \wedge \xi^1,
\]
\[\fl
d\xi^3 =
 \case{1}{4}  \, (\eta_1 +3\,\eta_2) \wedge \xi^3
+ 2\,\eta_3 \wedge \xi^1,
\]
\[
\fl
d\sigma_{11} =
 (\eta_1 -2\,\eta_2) \wedge \sigma_{11}
-4\, \eta_3 \wedge \sigma_{13}
- (\eta_4 -\theta_2) \wedge \theta_0
+ \eta_6 \wedge \xi^2
+ \eta_7 \wedge \xi^3
+ \eta_8 \wedge \xi^1
\]
\[
\fl\qquad\quad
-5\, \theta_1 \wedge (\xi^2+\sigma_{22})
+ \theta_2 \wedge \sigma_{12}
- 2\,\theta_3 \wedge\sigma_{23},
\]
\[
\fl
d\sigma_{12} =
 \case{1}{2}  \, (\eta_1 -3\,\eta_2)\wedge \sigma_{12}
-2\, \eta_3 \wedge \sigma_{23}
+ \eta_4 \wedge \xi^2
+ \eta_5 \wedge \xi^3
+ \eta_6 \wedge \xi^1
-2\, \theta_2 \wedge (\xi^2+\sigma_{22}),
\]
\[
\fl
d\sigma_{13} =
 \case{1}{4}  \, (3\,\eta_1 -7 \,\eta_2) \wedge 3\,\sigma_{13}
-3\, \eta_3 \wedge \sigma_{12}
+ \eta_5 \wedge \xi^2
+ \eta_6 \wedge \xi^3
+ \eta_7 \wedge \xi^1
-3\, \theta_3 \wedge (\xi^2+\sigma_{22}),
\]
\[
\fl
d\sigma_{22} =
\eta_4 \wedge \xi^1
- \case{1}{2}  \, \eta_1 \wedge \xi^2
- \case{1}{2}  \,  \eta_2 \wedge (3\,\xi^2+2\,\sigma_{22})
- \eta_3 \wedge \xi^3,
\]
\[
\fl
d\sigma_{23} =
 \case{1}{4}  \, (\eta_1 -5\,\eta_2) \wedge \sigma_{23}
- \eta_3 \wedge (\xi^2+\sigma_{22})
+ (\eta_4 -\theta_2) \wedge \xi^3
+ \eta_5 \wedge \xi^1,
\]
\[
\fl
d\eta_1 =
\xi^1 \wedge (\xi^2+ \sigma_{22}),
\]
\[
\fl
d\eta_2 = -3\,\xi^1 \wedge (\xi^2+ \sigma_{22}),
\]
\[
\fl
d\eta_3 =
 \case{1}{4}  \, (\eta_1 -\eta_2) \wedge \eta_3
+ \sigma_{23} \wedge \xi^1
-\xi^3 \wedge (\xi^2+ \sigma_{22}),
\]
\[
\fl
d\eta_4 =
 \eta_9 \wedge \xi^1
+ \case{1}{2}  \, (\eta_2 \wedge (3\,\theta_2-4\,\eta_4)+\eta_1 \wedge \theta_2)
+ \xi^2 \wedge \sigma_{22}
+ \xi^3 \wedge \sigma_{23},
\]
\[
\fl
d\eta_5 =
 \eta_9 \wedge \xi^3
+ \eta_{10} \wedge \xi^1
+ \case{1}{4}  \,(\eta_1 -9\,\eta_2)\wedge \eta_5
-3\, \eta_3 \wedge (\eta_4-\theta_2)
+ \sigma_{12} \wedge \xi^3
\]
\[
\fl\qquad\quad
+ 3\,\sigma_{23} \wedge (\xi^2+\sigma_{22}),
\]
\[
\fl
d\eta_6 =
 \eta_9 \wedge \xi^2
+ \eta_{10} \wedge \xi^3
+ \eta_{11} \wedge \xi^1
+ \case{1}{2}  \, (\eta_1 -5\,\eta_2) \wedge \eta_6
-4\, \eta_3 \wedge \eta_5
+ \sigma_{12} \wedge (7\,\xi^2+6\,\sigma_{22}),
\]
\[
\fl
d\eta_7 =
 \eta_{10} \wedge \xi^2
+ \eta_{11} \wedge \xi^3
+ \eta_{12} \wedge \xi^1
+ \case{1}{4}  \, (3\,\eta_1 -11\,\eta_2) \wedge \eta_7
-5\, \eta_3 \wedge \eta_6
+3\, \eta_4 \wedge \theta_3
\]
\[
\fl\qquad\quad
- (\eta_5 -3\,\theta_3)\wedge \theta_2
+3\, \sigma_{12} \wedge \sigma_{23}
+9\, \sigma_{13} \wedge (\xi^2+\sigma_{22}),
\]
\[
\fl
d\eta_8 =
(\eta_1 - 3\,\eta_2)\wedge \eta_8
+6\, (\eta_4 \wedge \theta_1
-\eta_3 \wedge \eta_7
+\sigma_{13} \wedge \sigma_{23})
+2\, \eta_5 \wedge \theta_3
-2\, (\eta_6-3\,\theta_1) \wedge \theta_2
\]
\begin{equation}
\fl\qquad\quad
- (\eta_9 +\sigma_{12}) \wedge \theta_0
+ \eta_{11} \wedge \xi^2
+ \eta_{12} \wedge \xi^3
+ \eta_{13} \wedge \xi^1
+ 12\,\sigma_{11} \wedge (\xi^2+\sigma_{22}).
\label{pKhZ_se}
\end{equation}
We have the following Maurer--Cartan forms
\[
\fl
\xi^1 =\frac{dt}{a},
\]
\[
\fl
\xi^2 =a\,\left(
\frac{u_{xxy}^2-u_x u_{xxx}^2}{u_{xxx}}\,dt+u_{xxx} dx+u_{xxy} dy
\right),
\]
\[
\fl
\xi^3 =2\,\frac{u_{xxy}}{u_{xxx}^{1/2}}\,dt +u_{xxx}^{1/2}\,dy,
\]
\[
\fl
\theta_2 = a^2\,u_{xxx}\,(du_x - u_{tx} \,dt- u_{xx}\,dx - u_{xy}\,dy),
\]
\[
\fl
\theta_3 = a^3\,u_{xxx}^{1/2}\,(
u_{xxx}\,(du_y - u_{ty} \,dt- u_{xy}\,dx - (u_{tx}+u_x\,u_{xx})\,dy)
\]
\[
\fl\qquad\qquad\qquad
-u_{xxy}\,(du_x - u_{tx} \,dt- u_{xx}\,dx - u_{xy}\,dy)
),
\]
\[
\fl
\eta_1 = 3\,\frac{da}{a}+2 \,\frac{du_{xxx}}{u_{xxx}}-u_{xx}\,dt,
\]
\[
\fl
\eta_2 = -\frac{da}{a}+3 \,u_{xx}\,dt,
\]
\begin{equation}
\fl
\eta_3 = a\,
\left(
\frac{du_{xxy}}{u_{xxx}^{1/2}}
-\frac{u_{xxy}\,du_{xxx}}{u_{xxx}^{3/2}}
+\frac{(u_{xx}u_{xxy}+u_{xy}u_{xxx})\,dt}{u_{xxx}^{1/2}}
+u_{xx}u_{xxx}^{1/2}\,dy
\right),
\label{pKhZ_MCfs}
\end{equation}
where $a \neq 0$ is a parameter. We do not need explicit expressions for the other
{\sc mc} forms of this pseudo-group in what follows.

\subsection{The Boyer--Finley equation}
For the symmetry pseudo-group of equation (\ref{BF}) we have the following structure equa\-ti\-ons
\[
\fl
d\theta_0 =
 \theta_0 \wedge (\theta_3-\sigma_{33})
+ \xi^1 \wedge \theta_1
+ \xi^2 \wedge \theta_2
+ \xi^3 \wedge \theta_3,
\]
\[
\fl
d\theta_1 =
 \eta_1 \wedge \theta_1
+\xi^1  \wedge \sigma_{11}
+\xi^2  \wedge \sigma_{33}
+\xi^3  \wedge \sigma_{13},
\]
\[
\fl
d\theta_2 =
\theta_2 \wedge (\eta_1 +\theta_3+\xi^3)
+\xi^1  \wedge \sigma_{33}
+\xi^2  \wedge \sigma_{22}
+\xi^3  \wedge \sigma_{23},
\]
\[
\fl
d\theta_3 =
\xi^1  \wedge \sigma_{13}
+\xi^2  \wedge \sigma_{23}
+(\theta_3 + \sigma_{33})\wedge \xi^3,
\]
\[
\fl
d\xi^1 =
(\sigma_{33}-\theta_3- \eta_1)\wedge \xi^1,
\]
\[
\fl
d\xi^2 =
(\eta_1 + \sigma_{33}+\xi^3)\wedge \xi^2,
\]
\[
\fl
d\xi^3 =
(\sigma_{33}- \theta_3) \wedge \xi^3,
\]
\[
\fl
d\sigma_{11} =
(2\, \eta_1 + \theta_3  -\sigma_{33})\wedge  \sigma_{11}
+ \eta_2 \wedge \xi^3
+ \eta_3 \wedge \xi^1
- \sigma_{13} \wedge \xi^2
+ \theta_1 \wedge (\xi^2+\sigma_{13}),
\]
\[
\fl
d\sigma_{13} =
(\eta_1+\theta_3-\sigma_{33}) \wedge \sigma_{13}
+ \eta_2 \wedge \xi^1
+ (\theta_3 +2\,\sigma_{33} -\xi^3)\wedge \xi^2,
\]
\[
\fl
d\sigma_{22} =
\sigma_{22} \wedge (2\,\eta_1+\theta_3+2 \, \xi^3+\sigma_{33})
+ \eta_4 \wedge \xi^3
+ \eta_5 \wedge \xi^2
+ \theta_2 \wedge (\xi^1+\sigma_{23})
- \sigma_{23} \wedge \xi^1,
\]
\[
\fl
d\sigma_{23} =
\sigma_{23} \wedge (\eta_1+\xi^3+\sigma_{33})
+ \eta_4 \wedge \xi^2
+ (\theta_3 +2\,\sigma_{33}-\xi^3)\wedge \xi^1,
\]
\[
\fl
d\sigma_{33} =
\xi^1 \wedge \sigma_{13}
+\xi^2 \wedge \sigma_{23}
+\xi^3 \wedge (\sigma_{33}-\theta_3),
\]
\[
\fl
d\eta_1 =
(\sigma_{13} +\xi^2) \wedge \xi^1,
\]
\[
\fl
d\eta_2 =
 \eta_6 \wedge \xi^1
+2\, (\eta_1 +\theta_3 - \sigma_{33}) \wedge \eta_2,
\]
\[
\fl
d\eta_3 =
 \eta_6 \wedge \xi^3
+ \eta_7 \wedge \xi^1
+ (3\,\eta_1 +2 \, (\theta_3-\sigma_{33}))\wedge  \eta_3
- \eta_2 \wedge (\theta_1+\xi^2)
-3\,\sigma_{11} \wedge (\xi^2+\sigma_{13}),
\]
\[
\fl
d\eta_4 =
 \eta_8 \wedge \xi^2
-2\,( \eta_1+\xi^3+\sigma_{33})\wedge \eta_4,
\]
\[
\fl
d\eta_5 =
 \eta_8 \wedge \xi^3
+ \eta_9 \wedge \xi^2
-(3\, \eta_1+\theta_3+3 \, \xi^3+2 \, \sigma_{33}) \wedge \eta_5
- \eta_4 \wedge (\theta_2+\xi^1)
\]
\begin{equation}
\fl\qquad\quad
-3\, \sigma_{22} \wedge (\xi^1+\sigma_{23}).
\label{BF_se}
\end{equation}
In what follows we need explicit expressions for the following {\sc mc} forms only:
\[
\theta_1 = \frac{u_{yy}}{a}\,(du_x - e^{u_y} u_{yy}\,dt - u_{xx}\,dx - u_{xy}\,dy),
\]
\[
\theta_3 =du_y - u_{ty}\,dt - u_{xy}\,dx - u_{yy}\,dy,
\]
\[
\xi^1 = a\,dt,
\]
\[
\xi^2 = \frac{e^{u_y}\,u_{yy}^2}{a}\,dx,
\]
\[
\xi^3 =u_{yy}\,dy,
\]
\[
\sigma_{33} =\frac{du_{yy}}{u_{yy}}+\theta_3,
\]
\begin{equation}
\eta_1 =-\frac{da}{a}+\frac{u_{ty}}{a}\,\xi^1+\sigma_{33}-\theta_3.
\label{BF_MCfs}
\end{equation}

\section{Deformed cohomologies and coverings}
Now we consider infinite-dimensional Lie algebras $\mathfrak{g}_1$ and $\mathfrak{g}_2$ defined by normal pro\-lon\-ga\-ti\-ons,
\cite{Cartan1,Cartan4,Vasilieva1972,Stormark2000}, of systems (\ref{pKhZ_se}) and (\ref{BF_se}), respectively, and study their second deformed cohomologies. In both cases they appear to be non-zero. Non-trivial 2-cocycles define in\-te\-gra\-ble extensions, \cite{BryantGriffiths,Morozov2009}, of (\ref{pKhZ_se}) and (\ref{BF_se}). Solutions to the integrable extensions coincide with the Wahlquist--Estabrook forms of known coverings of equations
(\ref{pKhZ}) and (\ref{BF}).

\subsection{The potential Khoklov-Zabolotskaya equation}

It is simple to show that $H^1(\mathfrak{g}_1)$ is generated by 1-form
\[
\zeta=3\,\eta_1+\eta_2.
\]
Denote by $I$ the exterior ideal generated by 1-forms
$\theta_i$, $0 \le i \le 3$,
$\xi^j$, $1 \le j \le 3$,
$\sigma_{11}$, $\sigma_{12}$, $\sigma_{13}$, $\sigma_{22}$, $\sigma_{23}$, $\eta_k$, $1 \le k \le 8$, that is, the left hand sides of equations
(\ref{pKhZ_se}) contain differentials of the generators of $I$.  Then direct computations show that
\[
\mathrm{dim}\,\left(H^2_{\lambda\,\zeta}(\mathfrak{g}_1) \cap I\right)
=
\left\{
\begin{array}{lll}
0, & \hspace{10pt}& \lambda \neq -\case{1}{4},
\\
1, & \hspace{10pt}& \lambda = -\case{1}{4},
\end{array}
\right.
\]
and
\[
H^2_{-\frac{1}{4}\,\zeta}(\mathfrak{g}_1) \cap I=
\mathbb{R}\,[\Omega],
\]
where
\[
\Omega=
\case{1}{4}\,(\eta_1-\eta_2-4\,\eta_3) \wedge (\xi^1+\xi^2+\xi^3)
+\theta_2 \wedge (\xi^1+\xi^3)
+\theta_3 \wedge \xi^1.
\]
We have a conjecture that
\[
\mathrm{dim}\,H^2_{\lambda\,\zeta}(\mathfrak{g}_1)
=
\left\{
\begin{array}{lll}
0, & \hspace{10pt}& \lambda \neq -\case{1}{4},
\\
1, & \hspace{10pt}& \lambda = -\case{1}{4}.
\end{array}
\right.
\]

From (\ref{d_deformed_2}) it follows that  equation
\[
d\omega -\case{1}{4}\,\zeta \wedge \omega = \Omega
\]
is compatible with system
(\ref{pKhZ_se}), that is, defines an integrable extension, \cite{BryantGriffiths,Morozov2009}, of  (\ref{pKhZ_se}).  Therefore
Lie's third inverse fundamental theorem in Cartan's form, \cite{Cartan1,Cartan4,Vasilieva1972,Stormark2000}, ensures existence of a solution $\omega \not \in I$ to the equation
\[
d\omega = \case{1}{4}\,(3\,\eta_1+\eta_2) \wedge  \omega
+
\]
\[
\qquad\quad
\case{1}{4}\,(\eta_1-\eta_2-4\,\eta_3) \wedge (\xi^1+\xi^2+\xi^3)
+\theta_2 \wedge (\xi^1+\xi^3)
+\theta_3 \wedge \xi^1.
\]
Since forms (\ref{pKhZ_MCfs}) are known, we can find $\omega$ explicitly:
\[
\omega = a^2\,u_{xxx}^{3/2}\,
\left(
d q - \left(\case{1}{3}\,q_x^3 - u_x\,q_x - u_y\right)\,dt -q_x\,dx - \left(\case{1}{2}\,q_x^2 - u_x\right)\,dy
\right).
\]
In this expression $q$ is a new variable (an ``integration constant''), while the new parameter $q_x$ can be expressed in terms of the free parameter $a$  of the {\sc mc} forms of the symmetry pseudo-group from the relation
\[
a=\frac{u_{xxx}^{1/2}}{u_{xxx}\,q_x - u_{xxy}}.
\]
The condition $\omega=0$ gives the covering system
\[
\left\{
\begin{array}{lll}
q_t &=& \case{1}{3}\,q_x^3 - u_x\,q_x - u_y,  \phantom{\frac{\frac{A}{A}}{\frac{A}{A}}}
\\
q_y &=& \case{1}{2}\,q_x^2 - u_x                 \phantom{\frac{\frac{A}{A}}{\frac{A}{A}}}
\end{array}
\right.
\]
for equation (\ref{pKhZ}).  This covering was found in \cite{Kuzmina} and then in
\cite{KupershmidtManin,Gibbons1985}.

\subsection{The Boyer-Finley equation}

For the symmetry pseudo-group of the Boyer-Finley equation we have  $H^1(\mathfrak{g}_2)=\mathbb{R}\,\zeta$ with
\[
\zeta=\sigma_{33}-\theta_3
\]
and
\[
\mathrm{dim}\,\left(H^2_{\lambda\,\zeta}(\mathfrak{g}_2) \cap I\right)
=
\left\{
\begin{array}{lll}
0, & \hspace{10pt}& \lambda \neq -1,
\\
2, & \hspace{10pt}& \lambda = -1,
\end{array}
\right.
\]
where $I$ is the exterior ideal generated by
$\theta_i$, $0 \le i \le 3$,
$\xi^j$, $1 \le j \le 3$,
$\sigma_{11}$, $\sigma_{13}$, $\sigma_{22}$, $\sigma_{23}$, $\sigma_{33}$, $\eta_k$, $1 \le k \le 5$, while
\[
H^2_{-\zeta}(\mathfrak{g}_2) =
\mathbb{R}\,[\Omega_1]
\oplus
\mathbb{R}\,[\Omega_2]
\]
with
\[
\Omega_1 = \eta_1\wedge (\xi^1+\xi^2+\xi^3)
-(\theta_1+\xi^2) \wedge \xi^1
+(\theta_3+\xi^3) \wedge \xi^2
\]
and
\[
\Omega_2= (\sigma_{33}- \theta_3) \wedge \xi^3.
\]
We have a conjecture that
\[
\mathrm{dim}\,H^2_{\lambda\,\zeta}(\mathfrak{g}_2)
=
\left\{
\begin{array}{lll}
0, & \hspace{10pt}& \lambda \neq -1,
\\
2, & \hspace{10pt}& \lambda = -1.
\end{array}
\right.
\]
The integrable extension which corresponds to $\Omega_1$
\[
d\omega = (\sigma_{33}-\theta_3) \wedge  \omega
+
\eta_1\wedge (\xi^1+\xi^2+\xi^3)
-(\theta_1+\xi^2) \wedge \xi^1
+(\theta_3+\xi^3) \wedge \xi^2
\]
has the following solution
\[
\omega =u_{yy}\,
\left(
dq-\left(u_t+e^{q_y}\right)\,dt+e^{u_y-q_y}\,dx-q_y\,dy
\right)
\]
with the relation $a =u_{yy}\,e^{q_y}$ between the free parameter $a$ and the new parameter $q_x$.
The corresponding covering
\[
\left\{
\begin{array}{lll}
q_t &=& u_t+e^{q_y},    \phantom{\frac{\frac{A}{A}}{\frac{A}{A}}}
\\
q_x &=& -e^{u_y-q_y}   \phantom{\frac{\frac{A}{A}}{\frac{A}{A}}}
\end{array}
\right.
\]
was found independently in \cite{Zakharov82,SavelievVershik,MalykhNutkuSheftelWinternitz}.

The solution to equation $d\omega = (\sigma_{33}-\theta_3) \wedge  \omega+ \Omega_2$ reads
$\omega = u_{yy} \,(dq + \ln u_{yy} \,dy)$ and is not interesting.

\section{Conclusion}
We have shown that for some {\sc pde}s their coverings arise quite naturally from the second deformed cohomologies of their symmetry pseudo-groups. It would be interesting to find out whether coverings of other {\sc pde}s can be derived using this construction. The further research will include clarification of relations between this technique and the other approaches to finding differential coverings.
This also leads to the question of other applications of deformed cohomologies of infinite-dimensional Lie algebras as well as to the problem
of improving the methods of their study.

\section{Acknowledgments}
I am very grateful to Boris Kruglikov for important and stimulating discussions at the initial stage of this work.

\end{document}